\documentclass[reqno,11pt]{amsart}
\usepackage{amsmath,amsfonts,amssymb,amsxtra,latexsym,amscd,enumerate,amsthm,verbatim}

\usepackage{hyperref}
\usepackage{graphicx}
\usepackage{physics}
\usepackage{bbm}

\usepackage{hyperref} 
\hypersetup{
    colorlinks=true,
    citecolor= blue,
    linkcolor=black,
    filecolor=magenta,      
    urlcolor=cyan
    }

\usepackage[margin=1in]{geometry}
\numberwithin{equation}{section}

\newcommand{\R}{\mathbb{R}}

\newcommand{\T}{\mathbb{T}}

\newcommand{\Z}{\mathbb{Z}}

\numberwithin{equation}{section} 

\newtheorem{theorem}{Theorem}[section]

\newtheorem{proposition}[theorem]{Proposition}
\newtheorem{corollary}[theorem]{Corollary}

\newtheorem{assumption}[theorem]{Assumption}

\usepackage[shortlabels]{enumitem}

\begin{document}

\title{Remarks on the rate of linear vortex symmetrization}

\author{Hao Jia}
\address{University of Minnesota}

\email{jia@umn.edu}


\thanks{HJ is supported in part by NSF grant DMS-1945179. }

\maketitle

\begin{abstract}
{\small}
We reformulate results from \cite{IonescuJia} in simplified forms and derive rigorously the bounds given in Bassom and Gilbert \cite{Ba1}, which provided interesting insights on the vortex symmetrization phenomenon. 

\end{abstract}

\section{Introduction and main results}
The dynamics of two dimensional incompressible fluid flows is often dominated by the interaction of vortices, through  symmetrization of small perturbations around a vortex and merging of same sign vortices. These processes are believed to play a significant role in reverse cascade of energy from small to large scales in two dimensional flows, which is often observed in numerical simulations. There is an extensive literature studying the dynamics of vortex solutions to two dimensional Euler and Navier Stokes equations with small viscosity, from both applied math and physics perspectives. See for example \cite{Ba1,Benzi,Brachet,McW1}. 

To state our motivations and results precisely, consider a steady state (vortex) given by a radial vorticity field $\Omega(r)$, with $\Omega'(r)>0$ for $r>0$ using the polar coordinate $(r,\theta)$. 
Let $U(r)e_{\theta}$ be the associated velocity field. Then from Biot-Savart law, $U(r)$ satisfies
\begin{equation}\label{Ur}
\partial_rU(r)+U(r)/r=\Omega(r),\,\, r\in(0,\infty), \quad{\rm with}\,\,\lim_{r\to0}U(r)=\lim_{r\to\infty}U(r)=0.
\end{equation}
The linearized two dimensional Euler equation around $U(r)e_{\theta}$ is
\begin{equation}\label{B1}
\partial_t\omega+\frac{U(r)}{r}\partial_{\theta}\omega-\frac{\Omega'(r)}{r}\partial_{\theta}\psi=0,\qquad (\theta,r)\in\mathbb{T}\times(0,\infty),
\end{equation}
with initial data $\omega_0(\theta,r)$. The stream function $\psi$ is determined through
\begin{equation}\label{B2}
\left(\partial_{rr}+r^{-1}\partial_r+r^{-2}\partial_{\theta}^2\right) \psi=\omega,\qquad (\theta,r)\in\mathbb{T}\times(0,\infty).
\end{equation}
We assume the orthogonality conditions
\begin{equation}\label{B3.1}
\int_{\R^2}\omega_0(\theta,r)(\cos{\theta},\sin{\theta})\,r^2drd\theta=0\qquad {\rm and }\qquad \int_{\mathbb{T}}\omega_0(\theta,r)\,d\theta=0.
\end{equation}
The general case can be reduced to the case with the orthogonality conditions \eqref{B3.1} by re-centering a slightly modified $\Omega(r)$.

Taking Fourier transform in $\theta$ in the equation \eqref{B1}-\eqref{B2} for $\omega$, we obtain that
 \begin{equation}\label{B4}
 \partial_t\omega_k+ik\frac{U(r)}{r}\omega_k-ik\frac{\Omega'(r)}{r}\psi_k=0, \quad (\partial_r^2+r^{-1}\partial_r-k^2r^{-2})\psi_k=\omega_k,
 \end{equation}
 for $k\in\mathbb{Z}, t\ge0, r\in(0,\infty)$, with initial data $\omega_0^k$. Here $\omega_k$ and $\psi_k$ are the $k-$th Fourier coefficients of $\omega$ and $\psi$. Below we assume that $k\in[1,\infty)$ unless otherwise specified, since the zero mode does not evolve and the case $k\in(-\infty,-1]$ can be studied similarly using a complex conjugation. 

Inspired by recent breakthroughs on asymptotic stability near shear flows for $2d$ Euler equations, see \cite{BM} on Couette flow and more recently \cite{IOJI_ActaMath,NaZh_Preprint} on general stable monotonic shear flow, precise results on the linearized dynamics around vortex solutions have been obtained in \cite{Bedrossian2019,IonescuJia,Ren2023}. See also \cite{Ionescu2022} for a result on nonlinear vortex symmetrization near a point vortex. These results show that the velocity and stream fields decay over time with explicit rates, and the vorticity field weakly converges to a radial function, a phenomenon often termed ``vortex symmetrization". One can also view the linearized asymptotic results as an essential step towards a proof of full nonlinear vortex symmetrization, a major open problem in 2d Euler equations. 

Linearized stability analysis has been studied, however, much earlier by Bassom and Gilbert \cite{Ba1}. It appears that, despite recent progresses, the results obtained in \cite{Ba1}, using matched asymptotics method and further verified by numerical computations, have not been fully proved.

There are two key predictions in \cite{Ba1}, which we reformulate slightly as follows.
\begin{enumerate}
    
    \item We can decompose the vorticity field for $r>0, t>1$ as
\begin{equation}\label{Pre1}
    \omega_k(t,r)=f_{0k}(t,r)e^{-ikU(r)t/r}+f_{1k}(t,r),
\end{equation}
where $f_{0k}(t,r)$ vanishes as $r\to 0+$ to the order $r^{\sqrt{k^2+8}}$ and $f_{1k}(t,r)$ decays over time with the rate $r^{-2-\sqrt{k^2+8}}t^{-1-\sqrt{k^2+8}}$. See (5.14a) in \cite{Ba1}.

\item For any $\eta\in C_0^\infty(\R^2)$, denoting $\eta_k(r)$ as the $k-$th mode of $\eta$ in $\theta$, we have for $t>1$,
\begin{equation}\label{Pre2}
    \int_{0}^{\infty}\omega_k(t,r)\eta_k(r)\,rdr\lesssim t^{-1-\sqrt{k^2+8}}.
\end{equation}
See figure 5 in \cite{Ba1} and the discussion below it. 
\end{enumerate}

The predicted decay rates over time in \eqref{Pre1}-\eqref{Pre2} are rather strong, while the existing results in \cite{Bedrossian2019,IonescuJia} only provide rigorous decay rates given by $1/t$. 

It turns out that one can derive \eqref{Pre1}-\eqref{Pre2} from the bounds on spectral density functions proved in \cite{IonescuJia} by Ionescu and the author. The paper \cite{IonescuJia} however focused its attention more on the optimal Gevrey regularity bounds, aiming at possible applications in nonlinear vortex symmetrization problem, and due to various technical reasons including the need to work with very precise Gevrey and weighted spaces, some of the estimates are quite involved. In this note we reformulate the main results in \cite{IonescuJia} in a more concise way (without insisting on the optimality of bounds at the Gevrey level), and rigorously justify the formulae \eqref{Pre1} and \eqref{Pre2} as a consequence of \cite{IonescuJia}.

To simplify our presentations as much as possible, we allow our implied constants to depend on the Fourier mode $k$. Our main assumptions on the background vorticity field $\Omega$ are the following. 

\begin{assumption}\label{MainAs}
We assume that the background radial vorticity profile $\Omega(r)$ satisfies the following natural conditions.
There exist constants $C_\ast\in(0,\infty)$ and $c_\ast\in(-\infty,0)$ such that for all $r\in(0,\infty)$ and $j\in\Z\cap[0,\infty)$ we have
\begin{equation}\label{Back}
\begin{split}
&0<\Omega(r)\leq \frac{C_\ast}{\langle r\rangle^6}, \quad \partial_r\Omega(r)<0,\quad\Big|(r\partial_r)^j \Big(\frac{\Omega'(r)}{r}-\frac{c_\ast}{\langle r\rangle^8}\Big)\Big|\leq \frac{C_\ast^{j+1}(j!)^{2}}{\langle r\rangle^{10}}r^2.
\end{split}
\end{equation}
In the above, we used the notation $\langle x\rangle:=\sqrt{x^2+2}$ for $x\in\R$, and the convention that $0!=1$.\\
\end{assumption}

We can now state our main result as follows.
\begin{theorem}\label{BMain}
Assume that $k\in\Z\cap[1,\infty)$, and that $\omega_{0k}(r)$ is smooth for $r\in(0,\infty)$ and satisfies the Gevrey regularity bounds in $v=\log r\in\R$, in \eqref{MPro3} and \eqref{MARr2}-\eqref{MARr3} below. Let $\omega_k(t,r)\in C^\infty([0,\infty)\times (0,\infty))$ be the unique solution to \eqref{B4} with initial data $\omega_{0k}$. Then we have the following conclusions. 
\begin{enumerate}
\item We can decompose the vorticity field for $r>0, t>1$ as
\begin{equation}\label{TPre1}
    \omega_k(t,r)=f_{k\rm loc}(t,r)e^{-ikU(r)t/r}+f_{k\rm nloc}(t,r),
\end{equation}
where $f_{k\rm loc}(t,r)$ satisfies the bounds for $r>0, t>1$ and $j\in[0,\infty)$,
\begin{equation}\label{TPre1.1}
    \big|(r\partial_r)^jf_{k\rm loc}(t,r)\big|\lesssim_{j,k} M_k\frac{r^{\sqrt{k^2+8}}}{1+r^{|k|+8+\sqrt{k^2+8}}},
\end{equation}
and $f_{k\rm nloc}(t,r)$ satisfies the bounds for $r>0, t>1$ and $j\in[0,\infty)$,
\begin{equation}\label{TPre1.2}
\begin{split}
 & \big|(r\partial_r)^jf_{k\rm nloc}(t,r)\big|\\
 &\lesssim_{j,k} M_k\frac{\big(r+r^{-1}\big)^2}{t}\min\bigg\{\big(r+r^{-1}\big)^{2\sqrt{k^2+8}}t^{-\sqrt{k^2+8}},\,\,1\bigg\}\frac{r^{\sqrt{k^2+8}}}{1+r^{|k|+\sqrt{k^2+8}+8}};
  \end{split}
\end{equation}

\item For any $\eta\in C_0^\infty(\R^2)$, with $\eta(r,\theta)=O(r^k)$ as $r\to 0$, we have for $t>1$,
\begin{equation}\label{TPre2}
    \Big|\int_{0}^{\infty}\omega_k(t,r)\eta_k(r)\,rdr\Big|\lesssim_{j,k} M_k\, t^{-1-\sqrt{k^2+8}}.
\end{equation}
 \end{enumerate}

 In the above, $M_k$ is the weighted Gevrey-2 norm of $\omega_{0k}$ in the variable $v=\log r$, given by \eqref{MARr2}-\eqref{MARr3}. The implied constants depend on the constants $C_\ast, c_\ast$ in \eqref{Back} and the constants from \eqref{MARr2}-\eqref{MARr3} in the definition of the Gevrey norms, as well as a structure constant associated with the linearized operator around $\Omega$, but are independent of $t>1, r\in(0,\infty)$. 
\end{theorem}

The assumptions \eqref{MARr2}-\eqref{MARr3} on the initial data are basically that $\omega_{0k}(r)$ is Gevrey-2 regular in $v=\log r$, that $\omega_{0k}(r)=\sigma_kr^k+O(r^{k+2})$ as $r\to0+$, and that $\omega_{0k}(r)=O(r^{-k-8})$ as $r\to+\infty$. 


Theorem \ref{BMain} provides a rigorous proof of the predictions \eqref{Pre1}-\eqref{Pre2} in \cite{Ba1}. Moreover, we have more precise information on the regularity of the ``profiles" $f_{k\rm loc}$ and $f_{k\rm nloc}$ of the vorticity field $\omega_k(t,r)$.

We note that one can also prove similar bounds to \eqref{TPre1.1}-\eqref{TPre1.2} with smooth but non-Gevrey initial data. Less smooth, such as Sobolev regular, initial data can also be considered with correspondingly weaker regularity control on the profiles. We use Gevrey regularity in our presentation since these results are already proved in \cite{IonescuJia} and can be used without additional proof. 

There is a similar vorticity depletion phenomenon for non-monotonic shear flows near critical points, investigated first in \cite{Bouchet}, and established rigorously in \cite{Wei2019} with sharp decay over time at the critical point for symmetric flows and qualitative bounds for general non-monotonic shear flows. A more quantitative version of the vorticity depletion for general non-monotonic shear flows was proved in \cite{Ionescu2023}

The rest of the paper is organized as follows. In section \ref{sec:old} we recall some general results on \eqref{B4} and essential estimates we need from \cite{IonescuJia}, and then we give the proof of Theorem \ref{BMain} in section \ref{sec:pmain}.
\section{Simplified bounds on the spectral density function}\label{sec:old}
 In this section, we introduce the spectral representation formula and recall the main bounds on the spectral density functions proved in \cite{IonescuJia} in simplified forms.


Define for $r\in\R^+:=(0,\infty)$,
\begin{equation}\label{B8}
b(r):=U(r)/r,\qquad d(r)=\Omega'(r)/r.
\end{equation}
It follows from \eqref{Ur}, \eqref{Back} and \eqref{B8} that for $r\in\R^+$,
\begin{equation}\label{Int0.1}
b(r)=U(r)/r=\int_0^1s\,\Omega(rs)\,ds\approx \frac{1}{\langle r\rangle^2},\qquad b'(r)\approx -\frac{r}{\langle r\rangle^4}, \qquad |b''(r)|\lesssim \frac{1}{\langle r\rangle^4}.
\end{equation}
In the above the implied constants depend only on the constants $c_\ast$ and $C_\ast$ in \eqref{Back}.

Define for $k\in\Z\backslash\{0\}, r,\rho\in (0,\infty)$, the function
\begin{equation}\label{B5}
G_k(r,\rho):=\left\{\begin{array}{ll}
                        \frac{\rho}{2|k|}\Big(\frac{r}{\rho}\Big)^{|k|}&\qquad{\rm for}\,\,r<\rho,\\
                         \frac{\rho}{2|k|}\Big(\frac{\rho}{r}\Big)^{|k|}&\qquad{\rm for}\,\,r\ge \rho.
 \end{array}\right.
\end{equation}
$G_k(r,\rho)$ is the Green's function for the differential operator $-\partial_r^2-r^{-1}\partial_r+k^2r^{-2}$ on $(0,\infty)$ with vanishing boundary conditions.

 For each $k\in\mathbb{Z}\backslash\{0\}$, we set for any $f\in X:=L^2\big(\R^+,rdr\big)$,
 \begin{equation}\label{B6}
 L_kf(r):=\frac{U(r)}{r}f(r)+\frac{\Omega'(r)}{r}\int_0^{\infty}G_k(r,\rho)f(\rho)\,d\rho.
 \end{equation}
  The equation \eqref{B4} can be reformulated as
 \begin{equation}\label{B7}
 \partial_t\omega_k+ikL_k\omega_k=0.
 \end{equation}
 
  Since $L_k$ is a compact perturbation of the simple multiplication operator $f\to b(r)f$, by general spectral theory, we can conclude that the spectrum of $L_k$ consists of the continuous spectrum $[0, b(0)]$ and possibly some discrete eigenvalues in $\R$. Our main assumptions \eqref{Back}-\eqref{B3} imply that there are no discrete eigenvalues. 
 
 

By standard theory of spectral projections, we then have
 \begin{equation}\label{B9}
 \begin{split}
 &\omega_k(t,r)=\frac{1}{2\pi i}\lim_{\epsilon\to0+}\int_{\R}e^{-ik\lambda t}\Big\{\left[(\lambda-i\epsilon-L_k)^{-1}-(\lambda+i\epsilon-L_k)^{-1}\right]\omega^k_0\Big\}(r)\,d\lambda\\
 &=\frac{-1}{2\pi i}\lim_{\epsilon\to0+}\int_{0}^{\infty}e^{-ikb(r_0) t}b'(r_0)\Big\{\big[(-b(r_0)-i\epsilon+L_k)^{-1}-(-b(r_0)+i\epsilon+L_k)^{-1}\big]\omega^k_0\Big\}(r)\,dr_0.
 \end{split}
 \end{equation}
We then obtain from \eqref{B4} that
 \begin{equation}\label{B10}
 \begin{split}
 \psi_k(t,r)&=\frac{1}{2\pi i}\lim_{\epsilon\to0+}\int_{0}^{\infty}e^{-ikb(r_0) t}b'(r_0)\int_0^{\infty}G_k(r,\rho)\\
 &\qquad\qquad\times\bigg\{\Big[(-b(r_0)-i\epsilon+L_k)^{-1}-(-b(r_0)+i\epsilon+L_k)^{-1}\Big]\omega^k_0\bigg\}(\rho)\,d\rho\, dr_0\\
 &=\frac{1}{2\pi i}\lim_{\epsilon\to0+}\int_{0}^{\infty}e^{-ikb(r_0) t}b'(r_0)\left[\psi_{k,\epsilon}^{-}(r,r_0)-\psi_{k,\epsilon}^{+}(r,r_0)\right]dr_0.
 \end{split}
 \end{equation}
 In the above, 
 \begin{equation}\label{B11}
 \begin{split}
 &\psi_{k,\epsilon}^{+}(r,r_0):=\int_0^{\infty}G_k(r,\rho)\Big[(-b(r_0)+i\epsilon+L_k)^{-1}\omega^k_0\Big](\rho)\,d\rho,\\
 &\psi_{k,\epsilon}^{-}(r,r_0):=\int_0^{\infty}G_k(r,\rho)\Big[(-b(r_0)-i\epsilon+L_k)^{-1}\omega^k_0\Big](\rho)\,d\rho.
 \end{split}
 \end{equation}
 We note that $\psi_{k,\epsilon}^{+}(r,r_0), \psi_{k,\epsilon}^{-}(r,r_0)$ satisfy for $\iota\in\{+,-\}$ and $r,r_0\in\R^+$,
 \begin{equation}\label{F7}
 \Big[k^2/r^2-r^{-1}\partial_r-\partial_r^2\Big]\psi_{k,\epsilon}^{\iota}(r,r_0)+\frac{d(r)}{b(r)-b(r_0)+i\iota\epsilon}\psi_{k,\epsilon}^{\iota}(r,r_0)=\frac{\omega^k_0(r)}{b(r)-b(r_0)+i\iota\epsilon}.
 \end{equation}

 It is more convenient to work in the variable $v$ with $r=e^v$ for $r\in(0,\infty)$. We therefore introduce for $k\in\mathbb{Z}\backslash\{0\}, \epsilon\in[-1/4,1/4]\backslash\{0\}$ and $\iota\in\{+,-\}$,
 \begin{equation}\label{B12}
\Pi^{\iota}_{k,\epsilon}(v,w):= \psi_{k,\epsilon}^{\iota}(r,r_0),\qquad{\rm where}\,\,r=e^v,\,\,r_0=e^w\,\,{\rm and}\,\,r,r_0\in(0,\infty),
 \end{equation}
 \begin{equation}\label{B12'}
 B(v):=b(r),\qquad D(v):=d(r),\qquad{\rm where}\,\,r=e^v,\,\,r\in(0,\infty).
 \end{equation}
 We also define 
 \begin{equation}\label{B12''}
 f_0^k(v):=\omega_0^k(r),\quad  \phi_k(t,v):=\psi_k(t,r), \qquad{\rm where}\,\,r=e^v,\,\,r\in(0,\infty).
 \end{equation}

 It follows from \eqref{B10}-\eqref{B12''} that for $\iota\in\{+,-\},v,w\in\R$, $k\in\Z\backslash\{0\}$, $\epsilon\in[-1/4,1/4]$,
 \begin{equation}\label{B13}
 k^2\Pi^{\iota}_{k,\epsilon}(v,w)-\partial_v^2\Pi^{\iota}_{k,\epsilon}(v,w)+e^{2v}D(v)\frac{\Pi^{\iota}_{k,\epsilon}(v,w)}{B(v)-B(w)+i\iota\epsilon}=e^{2v}\frac{f_0^k(v)}{B(v)-B(w)+i\iota\epsilon},
 \end{equation}
 and the normalized stream function has the representation formula
 \begin{equation}\label{B14}
 \phi_k(t,v)=\frac{1}{2\pi i}\lim_{\epsilon\to0+}\int_{\mathbb{R}}e^{-ikB(w)t}\partial_wB(w)\left[\Pi_{k,\epsilon}^{-}(v,w)-\Pi_{k,\epsilon}^+(v,w)\right]dw.
 \end{equation}
 The key is to study the regularity properties of the associated spectral density functions $\Pi_{k,\epsilon}^{-}(v,w)$ and $\Pi_{k,\epsilon}^{+}(v,w)$ for $v,w\in\R$, and $\epsilon\in[-1/4,1/4]\backslash\{0\}$ small.
 
 We summarize our calculations in the following proposition.
 \begin{proposition}\label{MPro} 
 Suppose $\Omega(r), r\in\R^+$, is a radial function satisfying the assumption \eqref{Back}. Let $U(r), r\in\R^+$, be given through
 \begin{equation}\label{MPro1}
\partial_rU(r)+U(r)/r=\Omega(r),\,\, r\in(0,\infty), \quad{\rm with}\,\,\lim_{r\to0}U(r)=\lim_{r\to\infty}U(r)=0.
 \end{equation}
 Consider $\Omega(r)$ as the steady vorticity profile for the two dimensional incompressible Euler equation with the associated velocity field $U(r)e_\theta$. The linearized equations around $\Omega$ for the perturbation vorticity $\omega(t,\theta,r)\in C^1([0,\infty), H^1(\T\times\R^+, r\langle r\rangle^8 d\theta dr))$ and the associated stream function $\psi(t,\theta,r)$, $t\ge0, \theta\in\mathbb{T}, r\in\R^+$ are given by
 \begin{equation}\label{MPro2}
 \begin{split}
& \partial_t\omega+\frac{U(r)}{r}\partial_{\theta}\omega-\frac{\Omega'(r)}{r}\partial_{\theta}\psi=0,\qquad\left(\partial_{rr}+r^{-1}\partial_r+r^{-2}\partial_{\theta}^2\right) \psi=\omega. \end{split}
\end{equation}
We assume that the initial vorticity deviation $\omega_0(\theta,r)$ satisfies the orthogonality conditions
\begin{equation}\label{B3}
\int_{\R^2}\omega_0(\theta,r)(\cos{\theta},\sin{\theta})\,r^2drd\theta=0,\qquad {\rm and \,\,for\,\,} r\in\R^+,\,\, \int_{\mathbb{T}}\omega_0(\theta,r)\,d\theta=0.
\end{equation}
For $k\in\Z\backslash\{0\}$ and $t\ge0, r\in\R^+$, letting 
\begin{equation}\label{MPro2.1}
\omega_k(t,r):=\frac{1}{2\pi}\int_{\T}\omega(t,\theta,r) e^{-ik\theta}\,d\theta,\quad \psi_k(t,r):=\frac{1}{2\pi}\int_{\T}\psi(t,\theta,r) e^{-ik\theta}\,d\theta, 
\end{equation} 
then $\omega_k(t,r), \psi_k(t,r)$ satisfy for $t\ge0, r\in\R^+$
 \begin{equation}\label{MPro2.5}
 \begin{split}
& \partial_t\omega_k+ik\frac{U(r)}{r}\omega_k-ik\frac{\Omega'(r)}{r}\psi_k=0,\quad\left(\partial_{rr}+r^{-1}\partial_r-k^2r^{-2}\right) \psi_k=\omega_k. \end{split}
\end{equation}
 Define for $t\ge0, r\in\R^+, v\in\R$ with $r=e^v$, the functions
 \begin{equation}\label{MPro3}
 \begin{split}
& b(r):=U(r)/r, \quad d(r):=\Omega'(r)/r, \quad B(v):=b(r),\quad D(v):=d(r),\\
& f_k(t,v):=\omega_k(t,r),\quad f_{0}^k(v):=\omega_{0}^k(r),\quad \phi_k(t,v):=\psi_k(t,r).
 \end{split}
 \end{equation}
 We have the representation formula for $t\ge0, v\in\R$,
 \begin{equation}\label{MPro4}
  \phi_k(t,v)=\frac{1}{2\pi i}\lim_{\epsilon\to0+}\int_{\mathbb{R}}e^{-ikB(w)t}\partial_wB(w)\left[\Pi_{k,\epsilon}^{-}(v,w)-\Pi_{k,\epsilon}^+(v,w)\right]dw,
 \end{equation}
 where the spectral density functions $\Pi_{k,\epsilon}^{\iota}(v,w)$ satisfy for $\iota\in\{+,-\}, v,w\in\R$, and $\epsilon\in[-1/8,1/8]\backslash\{0\}$ 
\begin{equation}\label{MPro5}
 k^2\Pi^{\iota}_{k,\epsilon}(v,w)-\partial_v^2\Pi^{\iota}_{k,\epsilon}(v,w)+e^{2v}D(v)\frac{\Pi^{\iota}_{k,\epsilon}(v,w)}{B(v)-B(w)+i\iota\epsilon}=\frac{e^{2v}f_0^k(v)}{B(v)-B(w)+i\iota\epsilon}.
 \end{equation}
 We also record the following bounds for later applications. There exists $C^{\ast}\in(0,\infty)$, depending on $C_\ast, c_\ast$ in \eqref{Back} such that for $v\in\R$ and $j\in\Z\cap[1,\infty)$, 
\begin{equation}\label{LAP100}
\begin{split}
&\Big|D(v)-\frac{c_\ast}{(1+e^{2v})^4}\Big|\leq C^\ast \frac{e^{2v}}{(1+e^{2v})^{5}}, \quad \big|\partial^j_vD(v)\big|\leq (C^\ast)^j(j!)^2 \frac{e^{2v}}{(1+e^{2v})^{5}},\\
&B(v)\approx \frac{1}{1+e^{2v}}, \quad \partial_vB(v)\approx -\frac{e^{2v}}{(1+e^{2v})^2}, \quad \big|\partial^j_vB(v)\big|\lesssim (C^\ast)^j(j!)^2\frac{e^{2v}}{(1+e^{2v})^2},\end{split}
\end{equation}
\begin{equation}\label{LAP101}
\partial_vB(v)=(c_\ast/4)e^{2v}+O\big(e^{4v}\big) \qquad{\rm for}\,\,v<0.
\end{equation}
 \end{proposition}

Denote $\alpha\wedge\beta:=\min\{\alpha,\beta\}$ for $\alpha, \beta\in\R$.  Define for $v, w\in\R$ the function $d(v,w)$ as
 \begin{equation}\label{MWF1In}
d(v,w):=\big|\big[\min\{v,w\}, \max\{v,w\}\big]\cap\big[\min\{w,0\}, 0\big]\big|.
\end{equation}

 
 Fix $\Phi_0(v)\in C^\infty(-\infty, -1)$ such that $\Phi_0\equiv 1$ on $(-\infty, -2]$.
 We choose also a smooth cutoff function $\Phi^\ast:\R\to[0,1]$ satisfying $\Phi^\ast\in C_0^{\infty}(-4,4)$ and $\Phi^\ast\equiv1$ on $[-2,2]$.

We make the following assumptions on the initial data. Recall that $c_\ast\in(-\infty,0)$ is from Assumption \ref{MainAs}.
 \begin{assumption}\label{MARr1}
There exist coefficients $\sigma_k\in\R$ for each $k\in\Z\backslash\{0\}$, and constants $M^\dagger_k\in(2\sigma_k,\infty), C_0>1$, such that the following statement holds. Define for $k\in\Z\backslash\{0\}$ and $v\in\R$,
 \begin{equation}\label{MARr2}
 F_{0k}(v):=f_0^k(v)-(\sigma_k/c_\ast) D(v)e^{|k|v}\Phi_0(v),
 \end{equation}
 then $F_{0k}$ satisfies the bounds for all $k\in\Z\backslash\{0\}$ and $v\in\R$,
 \begin{equation}\label{MARr3}
 \big|\partial_v^jF_{0k}(v)\big|\leq M_k  C_0^j(j!)^{2}\frac{e^{(|k|+2)v}}{1+e^{(2|k|+10)v}}.\\
 \end{equation}
\end{assumption}

Decompose for $k\in\Z\backslash\{0\}, \epsilon\in[-1/8,1/8]\backslash\{0\}$ and $v,w\in\R$,
\begin{equation}\label{MAR1}
\Pi_{k,\epsilon}^\iota(v,w):=(\sigma_k/c_\ast)e^{|k|v}\Phi_0(v)+\Gamma_{k,\epsilon}^{\iota}(v,w).
\end{equation}
It follows from \eqref{MPro5} that $\Gamma_{k,\epsilon}^{\iota}(v,w)$ satisfies the equation for $v,w\in\R, \epsilon\in[-1/8, 1/8]\backslash\{0\},$ $\iota\in\{+,-\}$, 
\begin{equation}\label{MAR2}
\begin{split}
&(k^2-\partial_v^2)\Gamma_{k,\epsilon}^{\iota}(v,w)+\frac{e^{2v}D(v)}{B(v)-B(w)+i\iota\epsilon}\Gamma_{k,\epsilon}^{\iota}(v,w)\\
&=\frac{e^{2v}F_{0k}(v)}{B(v)-B(w)+i\iota \epsilon}+(\sigma_k/c_\ast)\big(2|k|e^{|k|v}\partial_v\Phi_{0k}(v)+e^{|k|v}\partial_v^2\Phi_{0k}(v)\big).
\end{split}
\end{equation}

The main results we need below are the following precise estimates on the spectral density functions.

\begin{theorem}\label{MTH2}
Assume that $k\in\Z\backslash\{0\}$, $f_k(t,v), \phi_k(t,v)$ for $t\ge0, v\in\R$ are as in Proposition \ref{MPro}, and that the assumption \ref{MARr1} holds. Recall the definition \eqref{MWF1In} for $d(v,w)$ for $v,w\in\R$. Then we have the following conclusions.
 The limiting spectral density function 
\begin{equation}\label{MTH2.1}
\Gamma_k(v,w):=(-i)\lim_{\epsilon\to0+}\big[\Gamma_{k,\epsilon}^{+}(v,w)-\Gamma_{k,\epsilon}^{-}(v,w)\big]=2\lim_{\epsilon\to0+}\Im \,\Gamma_{k,\epsilon}^{+}(v,w)
\end{equation}
exists, as limit of functions in $L^2_{\rm loc}(\R^2)$. Define for $v, w\in\R$, the ``profile" $\Theta_k(v,w)$ for $\Gamma_k(v,w)$,
\begin{equation}\label{MTH2.10}
\Theta_k(v,w):=\Gamma_k(v+w,w).
\end{equation}
Then $\Gamma_k(v,w)$ and $\Theta_k(v,w)$ satisfy the following properties.

{\it (i) Bounds when $v$ is away from $w$.} For $v\in\R$ with $|v-w|\ge1/2$ and $w\in\R$, $\Gamma(v,w)$ is smooth in $v, w$, and for any $j_1,j_2\in\Z\cap[0,\infty)$,
\begin{equation}\label{MTH2.3}
\begin{split}
&\big|\partial_v^{j_1}\partial_w^{j_2}\Gamma_k(v,w)\big|\lesssim_{j_1,j_2,k} M_ke^{-|k||w|-(\sqrt{k^2+8}-|k|)d(0,w)-4\max\{w,0\}}e^{-|k||w-v|-(\sqrt{k^2+8}-|k|)d(v,w)}.
\end{split}
\end{equation}

{\it (ii) Bounds when $v$ is close to $w$.} For $w\in\R$ and $j\in\Z\cap[0,\infty)$, 
 we have
\begin{equation}\label{MTH2.5}
\begin{split}
&\Big\|\partial^j_w\Theta_k(v,w)\Big\|_{H^1(v\in[-2,2])}\lesssim_{j,k} M_ke^{-|k||w|-(\sqrt{k^2+8}-|k|)d(0,w)-4\max\{w,0\}}.
\end{split}
\end{equation}

{\it (iii) Equation for $\Gamma_k(v,w)$.} In addition, $\Gamma_k(v,w)$ satisfies for $v,w\in\R$ the equation
\begin{equation}\label{MTH2.6}
(k^2-\partial_v^2)\Gamma_k(v,w)+{\rm P.V.} \frac{e^{2v}D(v)\Gamma_k(v,w)}{B(v)-B(w)}=-2\pi\frac{e^{2w}\big(D(w)\digamma_k(w)-F_{0k}(w)\big)}{\partial_wB(w)}\delta(v-w),
\end{equation}
where ${\rm P.V.}$ represents principal value and $\digamma_k\in C^\infty(\R)$ satisfies the bound for $w\in\R$ and $j\in\Z\cap[0,\infty)$, 
\begin{equation}\label{MTH2.7}
\big|\partial_w^j\digamma_k(w)\big|\lesssim_{j,k} M_ke^{-|k||w|-(\sqrt{k^2+8}-|k|)d(w,0)}.
\end{equation}

{\it (iv) Refined regularity property of $\Theta_k(v,w)$ in $v$.} For any smooth $\zeta:[-10,10]\to \R$ and $|\zeta'|\gtrsim1$ on $[-10,10]$, define
\begin{equation}\label{MTH2.8}
\mathcal{H}(v, \rho, w):=\Theta(\zeta(v+\rho)-\zeta(\rho), w),
\end{equation}
we have for $j_1,j_2\in[0,\infty), \rho\in[-10,10]$, $w\in\R$,
\begin{equation}\label{MTH2.9}
\begin{split}
&\big\|\partial_\rho^{j_1}\partial_w^{j_2}\mathcal{H}(v,\rho,w)\big\|_{H^1(v\in[-1,1])}\lesssim_{j_1,j_2,k,\zeta} M_ke^{-|k||w|-(\sqrt{k^2+8}-|k|)d(w,0)-4\max\{w,0\}}.
\end{split}
\end{equation}

{\it (v) The explicitly solvable case $k=1$.} Assume that $k=1$. For $v, w\in\R$, we have 
\begin{equation}\label{BSD25}
\begin{split}
\Gamma_k(v,w)&=2\pi \frac{B(v)-B(w)}{(\partial_wB(w))^2}e^{v+w}\mathbf{1}_{v<w}\bigg\{f_0^k(w)-e^{-w}\frac{D(w)}{\partial_wB(w)}\int_{-\infty}^wf_0^k(\rho)e^{3\rho}\,d\rho\bigg\}.
\end{split}
\end{equation}

{\it (vi) Representation formula for the stream function and vorticity function.} Finally, we have the representation formula 
\begin{equation}\label{MTH2.91}
\begin{split}
&\phi_k(t,v)=-\frac{1}{2\pi }\int_\R e^{-ikB(w)t}\Gamma_k(v,w) \partial_wB(w)\,dw, \qquad f_k(t,v)=-e^{-2v}(k^2-\partial_v^2)\phi_k(t,v).\\
 \end{split}
\end{equation}

\end{theorem}

Theorem \ref{MTH2} is a simplified form of Theorem 1.4 in \cite{IonescuJia} (see also Proposition 6.1 in \cite{IonescuJia} for the formula \eqref{BSD25}), with a stronger decay  as $w\to+\infty$. The stronger decay can be seen in the factor $e^{-4\max\{w,0\}}$ from \eqref{MTH2.3}, \eqref{MTH2.5} and \eqref{MTH2.9}. Although the stronger bounds were not stated explicitly in \cite{IonescuJia} where it is not needed to get the slower $1/t$ decay, the proof follows from the same arguments for the proof of Proposition 8.1. The decay rate of the spectral density function $\Gamma_k(v,w)$, as $w\to\infty$, depends naturally on the decay rate of the background vortex and the initial data, and may decay very fast so that it is not a serious concern for us. The order of vanishing as $w\to-\infty$, corresponding to $r\to0+$ which ultimately leads to ``vorticity depletion", is a much more delicate issue, due to the nonlocal term in the linearized equation \eqref{B4}, and is one of the main questions that is addressed here and in \cite{IonescuJia}, as well as in other works \cite{Bedrossian2019,Ren2023}. 

The bounds \eqref{MTH2.3} and the representation formula \eqref{MTH2.91} show that the main contribution to the vorticity and stream functions $\omega_k, \psi_k$ are given by the spectral density function when the spectral parameter $w$ is close to the physical space parameter $v$, as expected. 

We refer the detailed proof of Theorem \ref{MTH2} to \cite{IonescuJia}, and only highlight two important observations used in the proof. The first key insight in the proof is to notice that for $w<<-1$ and $v\in[w+2,-10]$, the potential term in equation \eqref{MAR2} satisfies 
\begin{equation}\label{MTH4.0}
    \frac{e^{2v}D(v)}{B(v)-B(w)+i\iota\epsilon}\approx 8.
\end{equation}
The long range effect of the potential results in the faster decay described \eqref{MTH2.5} than $e^{-|k||w|}$ as $w\to-\infty$. The second key insight is to notice that $\Gamma_{k,\epsilon}^{+}(v,w)$ and $\Gamma_{k,\epsilon}^{-}(v,w)$ have significant cancellations, resulting in an extra decay when $v$ is away from $w$, as can be seen from the equation \eqref{MTH2.6}, which roughly speaking shows that $\Gamma_k(v,w)$ is ``generated" by the source term $-2\pi\frac{e^{2w}\big(D(w)\digamma_k(w)-F_{0k}(w)\big)}{\partial_wB(w)}\delta(v-w)$. The decay rates for $F_{0k}(w), D(w), \partial_wB(w)$ and $\digamma_k(w)$ imply that for $w>1$,
$$\Big|\frac{e^{2w}\big(D(w)\digamma_k(w)-F_{0k}(w)\big)}{\partial_wB(w)}\Big|\lesssim e^{-(|k|+4)|w|},$$
which heuristically explains the extra decay factor $e^{-4\max\{w,0\}}$ that we included in Theorem \ref{MTH2} in various bounds, in comparison with Theorem 1.4 of \cite{IonescuJia}. 

\section{Proof of the main results}\label{sec:pmain}
 We are now ready to give the proof of Theorem \ref{BMain}, using Theorem \ref{MTH2}. Consider first the bounds \eqref{TPre1.1}-\eqref{TPre1.2}. It follows from the representation formula \eqref{MTH2.91} and equation \eqref{MTH2.6} that for $t>1, v\in\R$,
 \begin{equation}\label{PMR1}
     \begin{split}
         f_k(t,v)=&\frac{1}{2\pi }\int_\R e^{-ikB(w)t}{\rm P.V.} \frac{D(v)\Gamma_k(v,w)}{B(v)-B(w)} \partial_wB(w)\,dw\\
           &+\Big[D(v)\digamma_k(v)-F_{0k}(v)\Big]e^{-ikB(v)t}. 
     \end{split}
 \end{equation}
 We define for $t>1,v\in\R$,
 \begin{equation}\label{PMR2}
     \begin{split}
         F_{k\rm loc}(t,v):=&\frac{1}{2\pi }D(v)\int_\R e^{-ik(B(w)-B(v))t}{\rm P.V.} \frac{\Gamma_k(v,w)}{B(v)-B(w)}\Phi^\ast(v-w) \partial_wB(w)\,dw\\
           &+\Big[D(v)\digamma_k(v)-F_{0k}(v)\Big],
     \end{split}
 \end{equation}
 and 
 \begin{equation}\label{PMR3}
     \begin{split}
         F_{k\rm nloc}(t,v):=\frac{1}{2\pi }D(v)\int_\R e^{-ikB(w)t} \frac{\Gamma_k(v,w)}{B(v)-B(w)}\big(1-\Phi^\ast(v-w)\big) \partial_wB(w)\,dw. 
     \end{split}
 \end{equation}
 Recalling the relation $r=e^v, v\in\R$, we then set for $t>1, r>0$,
 \begin{equation}\label{PMR4}
     f_{k\rm loc}(t,r):=F_{k\rm loc}(t,v),\quad f_{k\rm nloc}(t,r):=F_{k\rm nloc}(t,v). 
 \end{equation}
It suffices to prove the for $v\in\R, j\in[0,\infty)$, we have the bounds for $t\ge1$,
\begin{equation}\label{PMR5}
    \big|\partial_v^jF_{k\rm loc}(t,v)\big|\lesssim_{j,k}M_k\frac{e^{\sqrt{k^2+8}\,v}}{1+e^{(|k|+8+\sqrt{k^2+8})v}},
\end{equation}
and
\begin{equation}\label{PMR6}
\begin{split}
    &\big|\partial_v^jF_{k\rm nloc}(t,v)\big|\lesssim_{j,k}M_k\,e^{2|v|}t^{-1}\min\Big\{e^{2\sqrt{k^2+8}|v|}t^{-\sqrt{k^2+8}},\,1\Big\}\frac{e^{\sqrt{k^2+8}\,v}}{1+e^{(|k|+8+\sqrt{k^2+8})v}}. 
    \end{split}
\end{equation}
For \eqref{PMR6}, using the bounds \eqref{MTH2.3}-\eqref{LAP100}, and the observation that for $j_1,j_2\in[0,\infty)$ and $|v-w|\ge1$,
\begin{equation}\label{PMR6.1}
    \Big|\partial_v^{j_1}\partial_w^{j_2}\frac{1}{B(v)-B(w)}\Big|\lesssim_{j_1,j_2}\Big|\frac{1}{B(v)-B(w)}\Big|,
\end{equation}
we see that the proof for $j\ge1$ is similar to the case $j=0$. A useful inequality we can use is the following. Suppose $h(v,w)$ is smooth and compactly supported on $(j,j+1)$ satisfying for $\alpha\ge0$
\begin{equation}\label{PMR6.101}
    \big|\partial_w^\alpha h(v,w)\big|\lesssim_\alpha e^{-2\mu |w|}
\end{equation}
for some $j\in\Z$, and $\mu>0$, then for $t\ge1$,
\begin{equation}\label{PMR6.2}
    \Big|\int_\R h(v,w)e^{-ikB(w)t}\,dw\Big|\lesssim_\mu \min\Big\{\frac{1}{(|k|t)^\mu}, e^{-2\mu|j|}\Big\}.
\end{equation}
\eqref{PMR6.2} can be proved by a simple integration by parts argument. The desired bound \eqref{PMR6} then follow from \eqref{MTH2.3} and \eqref{PMR6.2} by a lengthy but straightforward calculation. 

To establish the bound \eqref{PMR5}, it suffices to consider the first term on the right hand side of \eqref{PMR2}. We follow the argument as in the proof of Theorem 1.5 in \cite{IonescuJia}. Assume that $|v-j|\leq1$ for some $j\in\Z$, we define the function $\zeta:[-4,4]\to \R$ such that we have the following identities for $|\rho-j|\leq 5$,
\begin{equation}\label{PMR8}
    \nu:=\frac{B(\rho)-B(j)}{MB'(j)},\qquad \rho:=\zeta(\nu)+j.
\end{equation}
In the above $M>1$ is chosen sufficiently large so that $\nu\in[-5,5]$ when $\rho\in[j-5,j+5]$. $\zeta$ then satisfies the assumption in (iv) of Theorem \ref{MTH2}. For $|v-j|\leq1, |w-j|\leq5$, denote
\begin{equation}\label{PMR9}
    \nu_1:=\frac{B(v)-B(j)}{MB'(j)},\quad \nu_2:=\frac{B(w)-B(j)}{MB'(j)}. 
\end{equation}
Then we can compute 
\begin{equation}\label{PMR10}
    \begin{split}
&\int_\R e^{-ik(B(w)-B(v))t}\,{\rm P.V.} \frac{\Theta_k(v-w,w)}{B(v)-B(w)}\Phi^\ast(v-w) \partial_wB(w)\,dw\\
&=\int_\R e^{-ik(\nu_2-\nu_1)MB'(j)t}\,{\rm P.V.} \frac{\Theta_k(\zeta(\nu_1)-\zeta(\nu_2),\zeta(\nu_2)+j)}{\nu_1-\nu_2}\Phi^\ast(\zeta(\nu_1)-\zeta(\nu_2)) \,d\nu_2\\
&=\int_\R e^{-ik\rho MB'(j)t}\,{\rm P.V.} \frac{\Theta_k(\zeta(\nu_1)-\zeta(\nu_1+\rho),\zeta(\nu_1+\rho)+j)}{-\rho}\Phi^\ast(\zeta(\nu_1)-\zeta(\nu_1+\rho))\, d\rho.
    \end{split}
\end{equation}
The desired bounds \eqref{PMR5} then follow from the smoothness of $\nu(v)$ in $v$, \eqref{MTH2.9} and smoothness of $\zeta$. 

We now turn to the proof of \eqref{TPre2} using \eqref{TPre1.1}-\eqref{TPre1.2}. This is similar to the argument in \cite{Ba1}. We first note that for $\eta\in C^\infty(0,\infty)$ satisfying for all $j\in[0,\infty)$,
\begin{equation}\label{PMR11}
\eta\equiv0 \,\,{\rm for}\,\,r\gg1, \quad |(r\partial_r)^j\eta(r)|\lesssim_{j,k} r^{k+2},
\end{equation}
the bounds for $t\ge1$ and $k\ge2$,
\begin{equation}\label{PMR12}
    \Big|\int_\R \omega_k(t,r)\eta(r)\,rdr\Big|\lesssim_{j,k} \frac{1}{t^{1+\sqrt{k^2+8}}}
\end{equation}
follows from \eqref{TPre1}-\eqref{TPre1.2}. For $k=1$, estimating the integral 
$$\int_\R \omega_k(t,r)\eta(r)\,rdr$$
directly using \eqref{TPre1}-\eqref{TPre1.2} seems to lead to a logarithmic loss in $t$. Similar issues also appeared when proving the bounds \eqref{TPre1}-\eqref{TPre1.2} in \cite{IonescuJia}. Remarkably, as observed already in \cite{Ba1} and later also in \cite{Bedrossian2019}, for $k=1$ the linearized equation \eqref{B4} is explicitly solvable. Using the formula \eqref{BSD25} from \cite{IonescuJia}, the desired bounds \eqref{TPre2} can be proved from direct computations. 

By the asymptotic expansion of $\eta_k(r)$ as $r\to 0+$, with leading order given by $r^k$ and $r^{k+2}$, it suffices to show for $t\ge1$,
\begin{equation}\label{PMR13}
    \Big|\int_\R r^k\omega_k(t,r)\varrho(r)\,rdr\Big|\lesssim_k \frac{1}{t^{1+\sqrt{k^2+8}}},
\end{equation}
if $\varrho\in C_0^\infty[0,\infty)$ with $\varrho\equiv1$ for $r\leq 1/8$. Using
\begin{equation}\label{PMR14}
    \begin{split}
        \int_\R r^k\omega_k(t,r)\varrho(r)\,rdr&=\frac{1}{2\pi}\int_{\R^2}\omega(t,r,\theta) \,\,\overline{r^k e^{ik\theta}\varrho(r)} \,rdr d\theta\\
        &=\frac{1}{2\pi}\int_{\R^2}\psi(t,r,\theta) \Delta\,\Big[\,\overline{r^k e^{ik\theta}\varrho(r)}\,\Big] \,rdr d\theta.
    \end{split}
\end{equation}
Since $r^ke^{ik\theta}$ is harmonic on $\R^2$, the function 
\begin{equation}\label{PMR15}
    H_k(r,\theta):=\Delta\,\Big[\,\overline{r^k e^{ik\theta}\varrho(r)}\,\Big]
\end{equation}
is compactly supported away from $r=0$. \eqref{PMR13} follows from the bounds for $t\ge1$ and $\varphi\in C_0^\infty(\R)$,
\begin{equation}\label{PMR16}
    \Big|\int_\R\phi_k(t,v)\varphi(v)\,dv\Big|\lesssim_k \frac{1}{t^{1+\sqrt{k^2+8}}}.
\end{equation}
\eqref{PMR16} can now be proved by decomposing $\phi_k(t,v)$, similar to $f_k(t,v)$, to the local and nonlocal part, and using bounds similar to \eqref{TPre1.1}-\eqref{TPre2} for the local and nonlocal parts respectively. Theorem \ref{BMain} is then proved.

\bigskip
\begin{center}
    {\bf Acknowledgement}\\
    We thank Rajendra Beekie for valuable discussions on the paper \cite{Ba1}.
\end{center}
 
\bibliographystyle{plain} 
\bibliography{bibliography}

\begin{thebibliography}{10}

\bibitem{Ba1}
A.~P. Bassom and A.~D. Gilbert.
\newblock The spiral wind-up of vorticity in an inviscid planar vortex.
\newblock {\em J. Fluid Mech.}, 371:109--140, 1998.

\bibitem{BM}
J.~Bedrossian and N.~Masmoudi.
\newblock Inviscid damping and the asymptotic stability of planar shear flows in the 2d euler equations.
\newblock {\em Publ. Math. Inst. Hautes Etudes Sci.}, 122:195--300, 2015.

\bibitem{Bedrossian2019}
J.~Bedrossian, M.~Coti Zelati, and V.~Vicol.
\newblock Vortex axisymmetrization, inviscid damping, and vorticity depletion in the linearized 2d euler equations.
\newblock {\em Annals of PDE}, 5(4), 2019.

\bibitem{Benzi}
R.~Benzi, G.~Paladin, S.~Patarnello, P.~Santangelo, and A.~Vulpiani.
\newblock Intermittency and coherent structures in two-dimensional turbulence.
\newblock {\em J. Phys. A: Math. Gen.}, 19:3771--3784, 1986.

\bibitem{Bouchet}
F.~Bouchet and H.~Morita.
\newblock Large time behavior and asymptotic stability of the 2d euler and linearized euler equations.
\newblock {\em Phys. D}, 239:948--966, 2010.

\bibitem{Brachet}
M.~Brachet, M.~Meneguzzi, H.~Politano, and P.~Sulem.
\newblock The dynamics of freely decaying two-dimensional turbulence.
\newblock {\em J. Fluid Mech.}, 194:333--349, 1988.

\bibitem{Ionescu2022}
A.~Ionescu and H.~Jia.
\newblock Axi‐symmetrization near point vortex solutions for the 2d euler equation.
\newblock {\em Communications on Pure and Applied Mathematics}, 75(4):818--891, 2022.

\bibitem{Ionescu2023}
A.~D. Ionescu, S.~Iyer, and H.~Jia.
\newblock On the stability of shear flows in bounded channels, ii: Non-monotonic shear flows.
\newblock {\em Vietnam Journal of Mathematics}, 2023.

\bibitem{IonescuJia}
A.~D. Ionescu and H.~Jia.
\newblock Linear vortex symmetrization: The spectral density function.
\newblock {\em Arch. Rational Mech. Anal.}, 246:61--137, 2022.

\bibitem{IOJI_ActaMath}
A.~D. Ionescu and H.~Jia.
\newblock Nonlinear inviscid damping near monotonic shear flows.
\newblock {\em Acta Mathematica}, 230(2):321--399, 2023.

\bibitem{NaZh_Preprint}
N.~Masmoudi and W.~Zhao.
\newblock Nonlinear inviscid damping for a class of monotone shear flows in finite channel.
\newblock {\em arXiv:2001.08564}, 2020.

\bibitem{McW1}
J.~McWilliams.
\newblock The vortices of two-dimensional turbulence.
\newblock {\em J. Fluid Mech.}, 219:361--385, 1990.

\bibitem{Ren2023}
S.~Ren, L.~Wang, D.~Wei, and Z.~Zhang.
\newblock Linear inviscid damping and vortex axisymmetrization via the vector field method.
\newblock {\em Journal of Functional Analysis}, 285(1):109919, 2023.

\bibitem{Wei2019}
D.~Wei, Z.~Zhang, and W.~Zhao.
\newblock Linear inviscid damping and vorticity depletion for shear flows.
\newblock {\em Annals of PDE}, 5(3), 2019.

\end{thebibliography}

\end{document}